\documentclass[10pt,times,twocolumn]{article}
\pdfoutput=1

\usepackage[all]{xy}
\usepackage{color,url}
\usepackage{subfigure}
\usepackage{bm,amsmath}
\usepackage{marvosym}
\usepackage{graphicx}
\usepackage{wrapfig}
\usepackage{epsf}
\usepackage{times,mathptm}
\usepackage{url}
\usepackage{algorithmic}
\usepackage{algorithm}
\usepackage{graphicx}
\graphicspath{Figures/}
\newcommand{\pv}{\mathbf{p}}
\newcommand{\pvs}{\pv^s}
\newcommand{\pvns}{\pv^{ns}}

\newcommand{\sv}{\mathbf{s}}

\newcommand{\fv}{\mathbf{f}}

\newcommand{\Gr}{\mathcal{G}}
\newcommand{\G}{{\Gr_g}}
\newcommand{\Lo}{{\Gr_l}}
\newcommand{\R}{{\Gr_r}}
\newcommand{\So}{\Gr_s}
\newcommand{\V}{\mathcal{V}}
\newcommand{\E}{\mathcal{E}}
\newcommand{\N}{\mathrm{Neb}}
\newcommand{\phase}{\mathbf{\theta}}
\newcommand{\ind}{\mathbf{x}}
\newcommand{\Lap}{\mathbf{L}}
\newcommand{\M}{\mathrm{M}}
\newcommand{\Mbig}{\mathbf{M}}
\newcommand{\flim}{\bar{\fv}}
\newcommand{\pvlim}{\bar{\pv}}

\newcommand{\propv}{\mathbf{\alpha}}
\newcommand{\pvr}{\pv^r}
\newcommand{\pvnot}{\pv^0}

\newcommand{\cps}[1]{\mathrm{csp}\left({#1}\right)}
\newcommand{\cli}[1]{\mathrm{cli}\left({#1}\right)}
\newcommand{\cg}[1]{\mathrm{cg}\left({#1}\right)}
\newcommand{\costt}[1]{\mathrm{cost}\left({#1}\right)}
\newcommand{\popt}{{\pvs}^*}
\newcommand{\sopt}{{\sv}^*}
\newcommand{\sr}{{\sv}^r}
\newcommand{\abs}[1]{\left|{#1}\right|}
\newcommand{\sI}{\mathrm{sI}}

\date{}
\begin{document}

\title{\Large \bf Operations-Based Planning for Placement and Sizing of Energy Storage\\
in a Grid With a High Penetration of Renewables}

%{\color{red} !!!!! Authors are needed only for the final submission - hide before submitting}
\author{{\bf Krishnamurthy Dvijotham}$^{(1)}$, {\bf Scott Backhaus} $^{(2)}$ and {\bf Misha Chertkov} $^{(3)}$\\
Department of Computer Science and Engineering, University of Washington, Seattle, WA 98195 $^{(1)}$\\
Material Physics and Applications Division $^{(2)}$ and \\
Center for Nonlinear Studies and Theoretical Division, $^{(3)}$\\
LANL, Los Alamos, NM 87545, USA}
\maketitle

\begin{center}
{\bf \large Abstract}
\end{center}
{\it
As the penetration level of transmission-scale time-intermittent renewable generation resources increases, control of flexible resources will become important to mitigating the fluctuations due to these new renewable resources.  Flexible resources may include new or existing synchronous generators as well as new energy storage devices.  The addition of energy storage, if needed, should be done optimally to minimize the integration cost of renewable resources, however, optimal placement and sizing of energy storage is a difficult optimization problem.  The fidelity of such results may be questionable because optimal planning procedures typically do not consider the effect of the time dynamics of operations and controls.  Here, we use an optimal energy storage control algorithm to develop a heuristic procedure for energy storage placement and sizing.  We generate many instances of intermittent generation time profiles and allow the control algorithm access to unlimited amounts of storage, both energy and power, at all nodes.  Based on the activity of the storage at each node, we restrict the number of storage node in a staged procedure seeking the minimum number of storage nodes and total network storage that can still mitigate the effects of renewable fluctuations on network constraints.  The quality of the heuristic is explored by comparing our results to seemingly ``intuitive'' placements of storage. 

}
\vspace{0.5cm}

\section{Introduction}
\label{sec:Intro}
Electrical grid planning has traditionally taken two different forms; operational planning and expansion or upgrade planning.  The first is concerned with the relatively short time horizon of day-ahead unit commitment or hour-ahead or five-minute economic dispatch.  The focus is on controlling assets that are already present on the system to serve loads at minimum cost while operating the system securely.  The second typically looks out many years or decades and is focused on optimal addition of new assets, again with a focus on minimizing the cost of electricity.  When a system consists entirely of controllable generation and well-forecasted loads, the network power flows do not deviate significantly or rapidly from well-predicted patterns.  In this case, expansion planning can be reasonably well separated from operational planning.  In the simplest approach, expansion plans may be optimized against two extreme cases, e.g. the system's maximum and minimum load configurations.

As the penetration of time-intermittent renewables increases, expansion and operational planning will necessarily become more coupled.  For an electrical grid with large spatial extent, renewable generation fluctuations (here, we focus on wind generation) at well-separated sites will be uncorrelated on short time scales\cite{Gibescu2009,08LBNL}, and the intermittency of this new non-controllable generation will cause the patterns of power flow to change on much faster time scales than before, and in unpredictable ways.  New equipment (e.g. combustion turbines or energy storage) and control systems may have to be installed to mitigate the network effects of renewable generation fluctuations to maintain generation-load balance.  It is at this point where operations planning must interface with expansion planning.  The optimal placement and sizing of the new equipment, if required, depends on how the rest of the network and its controllable components respond to the fluctuations of the renewable generation.  Overall, we desire to install a minimum of new equipment by placing it at network nodes that afford us a high degree of controllability, i.e. nodes where controlled power injection and/or consumption have a significant impact on the network congestion introduced by the renewable fluctuations.  From the outset, it is not clear which nodes provide the best controllability.  Implicitly, placing a minimum of new equipment implies that it will experience a high duty, thus avoiding outcomes where equipment is only used for a small fraction of time.

Before discussing our initial approach at integrating operational planning and expansion planning, we first summarize a few methods for mitigating the intermittency of renewable generation.  When renewable penetration is relatively low and the additional net-load fluctuations are comparable to existing load fluctuations, a power system may continue to operate ``as usual'' with primary and secondary regulation reserves\cite{99HK} being controlled via a combination of distributed local control, i.e. frequency droop, and centralized control, i.e. automatic generation control (AGC).  In this case, {\it planning} for renewables may simply entail increasing the level of reserves to guard against the largest expected fluctuation in {\it aggregate renewable output}.

As the penetration level grows, system operators may simply continue to increase the level of reserves, however, if this fast-moving generation comprises natural gas combustion turbines, this simple planning will generally result in increased renewable integration costs\cite{Meibom2010} which are usually spread over the rate base.  Alternatively, operational planning can be improved by incorporating long-range and short-range renewable generation forecasting to better schedule the controllable generation (energy and reserves) to meet net load and operate reliably\cite{Meibom2010,Bouffard2008,BPA-wind-hirst}.  Long-range forecasts, typically day-ahead or farther, are used in a unit commitment optimization to ensure adequate generation and reserves will be online to meet the expected net load and its uncertainty.  Short-range forecasts, typically hour ahead or shorter, are used in an economic dispatch optimization that sets actual generation and reserve levels of the committed generation.  Some have investigated rolling unit commitment\cite{Tuohy2007,Meibom2010} where updated wind forecasts are used modify the unit commitment more frequently.  Simulations have resulted in lower overall renewable integration costs.

Both unit commitment and economic dispatch seek minimize the cost of electricity, however, they must also respect system constraints including generation limits, transmission line thermal limits, voltage limits, system stability constraints, and N-1 contingencies.  Previous works\cite{Meibom2010,Bouffard2008,Tuohy2007,BPA-wind-hirst} have generally looked at the effects of stochastic generation on the economics and adequacy of {\it aggregate} reserves while not considering such network constraints.  These constraints may be respected for a dispatch based on a {\it mean} renewable forecast, however, if the number of renewable generation sites and their contribution to the overall generation is significant, verifying the system security of all probable renewable fluctuations (and the response of the rest of the system) via enumeration will be exponentially complex in the number of sites resulting in a computationally intractable problem.

The approaches summarized above do not consider network constraints or the behavior of the system on time scales shorter than the time between economic dispatches (one hour in the case of \cite{Meibom2010}).  In particular, they do not model how fast changes in renewable generation and the compensating response of regulation reserves interact with network constraints.  In our initial study, we augment the approaches summarized above by focusing on the behavior of the electrical network at a finer time resolution and investigate how the control of energy storage affects its placement and sizing.  We presume that the unit commitment problem has been solved, and at the start of a time period, we dispatch the controllable generation by solving a DC optimal power flow  (DCOPF) based on the {\it site-specific mean forecast} for wind generation.  In the time before the next DCOPF is executed, the wind generation fluctuates, and we control a combination of existing generators (using a simplified description of frequency droop control) and energy storage to avoid violations of network constraints and generation limits.  For each level of wind penetration, we generate many different realizations of wind fluctuations, and we gather statistics on the activity of the energy storage at each node.  The statistics from simulated system operations are then coupled to the expansion planning process by developing a heuristic to guide the optimal placement and sizing of storage throughout the network--a result that cannot be achieved with the previous approaches described above.

The rest of this manuscript is organized as follows.  Section~\ref{sec:math} lays out the system model including the development of our heuristic for placement and sizing of energy storage.  Section~\ref{sec:exp} describes the simulations carried out on a slightly modified version of RTS-96\cite{RTS96} and the application of our design heuristic.  Section~\ref{sec:inter} gives some interpretations of the results, and Section~\ref{sec:con} wraps up with some conclusions and directions for future work.

\section{Mathematical Formulation}
\label{sec:math}
We briefly describe the mathematical formulation of the optimal control problem and the algorithm used to solve it. We begin with a quick review of the DC power flow model and the proportional control scheme used to emulate frequency droop and AGC\cite{Chertkov2011}.
\subsection{Background and Notation}
Let $\Gr=(\E,\V)$ denote the undirected graph underlying the power grid, and $n$ denote the number of nodes in the grid. Nodes are classified into three types:

\begin{itemize}
\item[1] Loads $\Lo$(assumed fixed over time scale of the study relative to the renewable generation).
\item[2] Traditional Generators $\G$ (whose output can be controlled).
\item[3] Renewable generators $\R$ (generators based on renewable sources whose output fluctuates over time).
\end{itemize}

Given the vector $\pv$ of power generated/consumed at every node in the grid, the power flow equations determine how much power $\fv_{ij}$ flows through each line in the network from node $i$ to node $j$. For any $S\subset \V$, we denote the power generation at nodes in $S$ by $\pv_S=\{\pv_i:i \in S\}$ and use similar notation for other quantities(like energy).The DC power flow equations are a linearization of the exact AC power flow equations, and have been widely used to compute optimal economic dispatch of generators in the grid.  If $\N(j)$ denotes the set of neighbors of node $j$ in the grid and $\ind_{ij}$ denote the inductance of the line from $i$ to $j$, the DC power flow equations are given by:
\begin{align}
 \pv_i & = \sum_{j \in \N(i)}\frac{\phase_j-\phase_i}{\ind_{ij}} \nonumber \\
 \fv_{ij} & = \ind_{ij}(\phase_j-\phase_i), \nonumber
\end{align}
where $\phase_i$ is the phase of the complex voltage at node $i$.

Written in matrix form, the first set of equations become $\pv=-\Lap_\ind\phase$ where $\Lap_\ind$ is the Laplacian of the graph $\Gr$ with edge weights $1/\ind_{ij}$. The Laplacian is not invertible, but if we are given a balanced power configuration, i.e., $\sum_i \pv_i=0$, we can always find a unique solution for the $\phase$. This mapping can be constructed from the eigenvalue decomposition of $\Lap_\ind=\sum_i w_i e_i e_i^T $ by inverting the non-zero eigenvalues to get $\M_\ind=\sum_{i:w_i>0} \frac{1}{w_i} e_ie_i^T$ and $\phase=\M_\ind\pv$. Letting ${\M_\ind}^i$ denote the $i$-th row of $\M_\ind$, we can summarize the power flow equations as:
\begin{align}
 \fv	=	\Mbig\pv \nonumber \\
 \Mbig^{ij}	=	\ind_{ij}({\M_\ind}^i-{\M_\ind}^j) \label{eq:DCPF}
\end{align}

The aim of our energy storage control scheme is to dispatch power from storage such that, in spite of fluctuations in the renewable generation, the following constraints are maintained all the time:
\begin{align}
 -\flim & \leq \Mbig\pv && \leq \flim && \text{(Transmission Capacities)}\nonumber \\
 0 & \leq \pv_\G && \leq \pvlim_\G && \text{(Generation Capacities)} \nonumber \\
 & \sum_i \pv_i &&=0 && \text{(Power Balance)} \label{eq:Contraint}
\end{align}
where the sum over $i$ in the third constraint is a sum over three types of nodes.

\subsection{Power Generation and Proportional Control}
We will design a control scheme that operates in the duration $T$ between two successive DCOPFs. At the beginning of each such time horizon, we assume that a subset of active generators $\G$ are committed (based on the predetermined unit-commitment).  We set their generation levels according to a DCOPF:
\begin{align*}
 \min_{\pv_\G} & \sum_{i \in \G} w_i \pv_i \\
 \textrm{Subject to Eqs.~\ref{eq:DCPF} and \ref{eq:Contraint}}
\end{align*}
where $w_i$ are generation costs for each generator. Let $\pvnot$ be the complete vector of initial generation at every node based on the DCOPF. The output $\pv_\R$ of renewabnle generation $\R$ is taken to be the mean of the renewable generation $\pvnot_\R$ over the time between DCOPFs.  We assume we know this value exactly, which implies we have perfect predictions of the mean output of the renewables.

In reality, the renewable generation $\pv_\R$ fluctuates during this time causing a mismatch between total generation and total load. We assume that the generators $\G$ respond by changing their output by dividing the mismatch in some predetermined proportion $\propv$ (usually chosen based on generation capacities)\cite{Chertkov2011}. For notational convenience, we take $\propv$ to be a vector of length $n$, with entries $\propv_i=0$ for $i \not\in\G$. Let $\pvr=\pv_\R-\pvnot_\R$ denote the vector of renewable generation change. We assume for simplicity that loads are fixed, which is reasonable since in many cases fluctuations in renewables will be much larger than fluctuations in loads. The generators $\G$ will then respond as:
\[ \pv_\G=\pvnot_\G+\propv_\G(-\sum \pvr), \]
where non-zero components of $\propv_\G$ associated with generators are chosen such that the net balance of power is maintained, $\sum \pvr=\sum_{i\in\G} \pv_i$.

\subsection{Controlling Storage and Constraints on Storage}

We use $\sv$ to denote the vector of energies stored at each node. However, the quantity we control directly is $\pvs$, the power drawn from each storage node at any given time.  In our storage placement and sizing heuristic, we allow $\pvs$ and $\sv$ to be as high as needed to ensure that the constraints \eqref{eq:Contraint} are satisfied. We also restrict the set of nodes that have storage to a set $\So$. However, again for notational convenience, $\pvs,\sv$ are still assumed to be vectors of length $n$ and we set $\pvs_{j}=0,\sv_{j}=0$ for $j \not\in \So$.

In the paper, we use the optimal energy storage control as a heuristic to decide placement and sizing of storage. This decision only assumes that we know what the renewable generation is going to be over the next window of time of length $T$ until the next DCOPF. Thus, we would like this decision to be robust to the behavior of wind beyond $T$ and that our control over the next time window can work independently of what the control was in the current time window. Hence we enforce the following invariance condition: No net energy is exchanged between the grid and the storage during the window of time $T$. Mathematically, this means that we require:
\begin{align}
 \sv(0)=\sv(T) \label{eq:Cyclic}
\end{align}
We also assume that the storage is not accounted for in the DCOPF, i.e $\pvs(0)=0$.

\subsection{Overall Power generation}
If storage is being dispatched at a node, it adds to the pre-existing generation $\pvns$ such that total generation is $\pv=\pvns+\pvs$, where $\pvns$ is the generation without including the local storage dispatch.  For $\Lo$ or $\R$, $\pvns$ is just the fixed load value or the time-variable renewable generation.  The storage dispatch simply adds to these outputs.

For generators, i.e. $\G$, the situation is slightly more complicated.  The dispatch of storage modifies the instantaneous load-generation imbalance so that the proportional control we proposed for $\G$ is slightly modified.  The controlled generation at a node responds to the sum of storage dispatch at all nodes, including its own, i.e.
\[\pvns_\G=\pvnot_\G+\propv_\G\left(-\sum_{i \in \R} \pvr_i -\sum_i \pvs_i\right).\]
Thus, the overall power generation from a $\G$ node is
\begin{equation}\label{eq:pgen-and-storage}
\pv = \pvnot+\pvs+\pvr+\propv\left(-\sum_{i \in \R} \pvr_i -\sum_i \pvs_i\right)
\end{equation}

\subsection{Optimal Control with perfect Forecasts}
\begin{figure}[htb]
\begin{center}
 \includegraphics[width=.4\textwidth,height=.15\textheight]{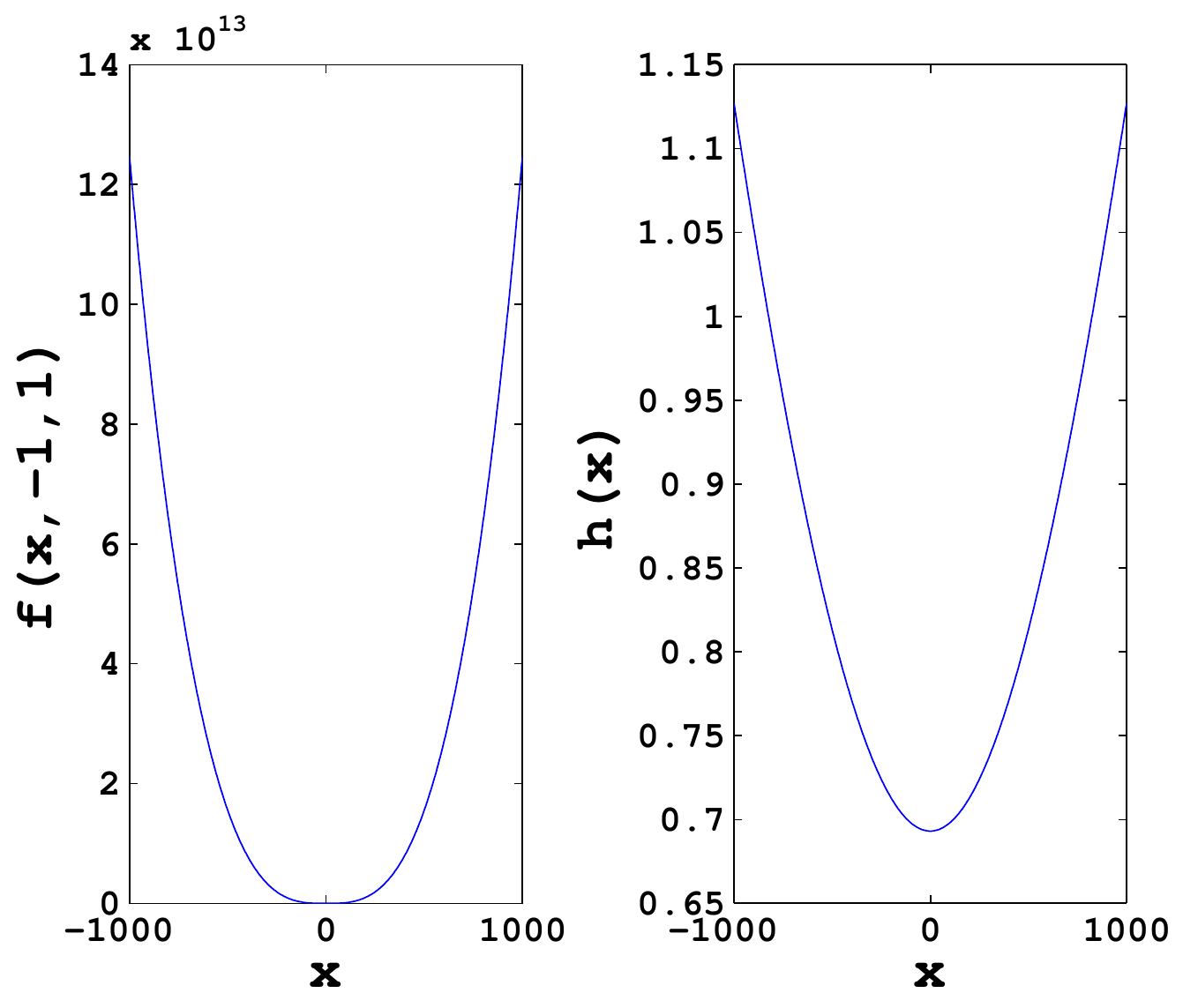}
\caption{Penalty Functions} \label{fig:Penalty}
\end{center}
\end{figure}

In this initial work on operations-based placement and sizing of energy storage, we implement an optimal control algorithm that assumes perfect forecasts of renewable generation, i.e. the output of each renewable generation node is known exactly over the time horizon $T$.  The control algorithm uses this information to decide how much power to draw from or inject into each storage node in order to keep the grid within the constraints \eqref{eq:Contraint}.  We note that, at the same time, the generation $\G$ is responding according to \eqref{eq:pgen-and-storage}. The control algorithm computes this optimal dispatch by minimizing a cost function consisting of three terms. The cost function tries to enforce the constraints \eqref{eq:Contraint} in a soft manner: We do this because in some situations, it might be infeasible to maintain all the constraints strictly. In those cases, we allow violations to occur but the controller suffers a cost that blows up cubically with the magnitude of the constraint violation. Thus, the controller will do its best to keep the system within its constraints, but when it is infeasible to do so, it will allow the violations to occur.

To make the notation compact, we define two auxillary functions $f$ and $h$, which are plotted in Figure \ref{fig:Penalty} and defined below.
\[
  f(x,a,b) = \left\{
  \begin{array}{l l}
    0 & \quad \text{if $a \leq x \leq b$}\\
    (\kappa_f*(x-b))^3 & \quad \text{if $x>b$}\\
    (\kappa_f*(x-a))^3 & \quad \text{if $x<a$}
  \end{array} \right.
\]
\begin{align*}
h(x) & = \log\left(\cosh\left(\kappa_h x\right)\right)
\end{align*}
where $\kappa_f=50,\kappa_h=.001$.

The three costs in the penalty function are:
\begin{itemize}
 \item[1] $\cli{\pv}$: penalizes overloading of transmission lines:\\
$\cli{\pv} = \sum_{(i,j) \in \E} f(\Mbig^{ij}\pv,-\flim_{ij},\flim_{ij})$
\item[2] $\cg{\pvns}$: penalizes violation of capacity constraints of controllable generators:\\
$\cg{\pvns} = \sum_{i \in \G} f(\pvns_i,0,\pvlim_i)$
\item[3] $\cps{\pvs}$: penalizes the absolute value of power drawn from storage (increases very slowly):\\
$\cps{\pvs} = \sum_{i} h(\pvs_i)$\\
The fact that $h$ increases very slowly ensures that $\cli{\pv},\cg{\pvns}$ dominate this cost by far, so that the controller will never choose to allow violations to occur in order to save on power capacity.
\end{itemize}
We define the total cost as $\costt{\pv,\pvns,\pvs}=\cli{\pv}+\cg{\pvns}+\cps{\pvs}$. The best control strategy minimizes this cost function integrated over time.  In order to leverage practical numerical techniques, we discretize the time axis uniformly with time step $\Delta$ to get $T_f=\frac{T}{\Delta}$ time steps and assume that our controls are constant over each time period $\Delta$.  We seek
\begin{align*}
 &\min_{\pvs_{\So}(0),\ldots,\pvs_{\So}(T_f-1)} \sum_{t=0}^{T_f} \cli{\pv}+\cg{\pvns}+\cps{\pvs} \\
 &\textrm{subject to}  \\
 &\sv(t) = \sum_{\tau=0}^{t-1} \pvs(t)\Delta \\
 &\sv(T_f)=\sv(0)   \\
 &\pvs_i(t) = 0 \quad i \not\in \So  \\
 &\pv(t) = \pvnot+\pvs(t)+\pvr(t)+\propv\left(-\sum_{i \in \R} \pvr_i(t) -\sum_i \pvs_i(t)\right) \\
 &\pvns(t) = \pv(t)-\pvs(t) \\
\end{align*}

Our problem falls within the well-studied class of deterministic, discrete-time optimal control problems and can be solved using standard nonlinear programming techniques \cite{betts2001practical}. Note that $\pvnot$ and $\pvr$ are known from the DCOPF and the perfect renewable forecasts, respectively. Thus, we can express $\pv(t),\pvs(t),\pvns(t)$ directly in terms of $\pvs_{\So}(t)$ using the second, third and fourth constraints above to get $\pv(\pvs_{\So}),\pvs(\pvs_{\So}),\pvns(\pvs_{\So})$. Further we can write $\pvs_{\So}(t)=\frac{\sv_{\So}(t+1)-\sv_{\So}(t)}{\Delta}$, with the convention that $\sv_{\So}(T_f)=\sv_{\So}(0)$ (in order to satisfy the constraint \eqref{eq:Cyclic}). This leaves us with the following unconstrained optimization problem:
\begin{align}
  \min_{\sv_{\So}(0),\ldots,\sv_{\So}({T_f}-1)} \sum_{t=0}^{T_f} &\cli{\pv\left(\frac{\sv_{\So}(t+1)-\sv_{\So}(t)}{\Delta}\right)}+\nonumber\\
    &\cg{\pvns\left(\frac{\sv_{\So}(t+1)-\sv_{\So}(t)}{\Delta}\right)}+\nonumber\\
	&\cps{\left(\frac{\sv_{\So}(t+1)-\sv_{\So}(t)}{\Delta}\right)}\label{eq:OptProb}
 \end{align}

It can be shown that this is a convex optimization problem \cite{boyd2004convex} that can be solved efficiently using Newton's method (since the Hessian has a sparse block structure arising from the pairwise interaction). In practice, a Levenberg-Marquardt correction \cite{nocedal1999numerical} is required to ensure convergence as the Hessian can become numerically singular. We use software from \cite{minFunc} to solve the problem .

\subsection{Optimal Sizing and Placement of Storage}
\begin{figure}[htb]
\begin{center}
 \includegraphics[scale=0.45]{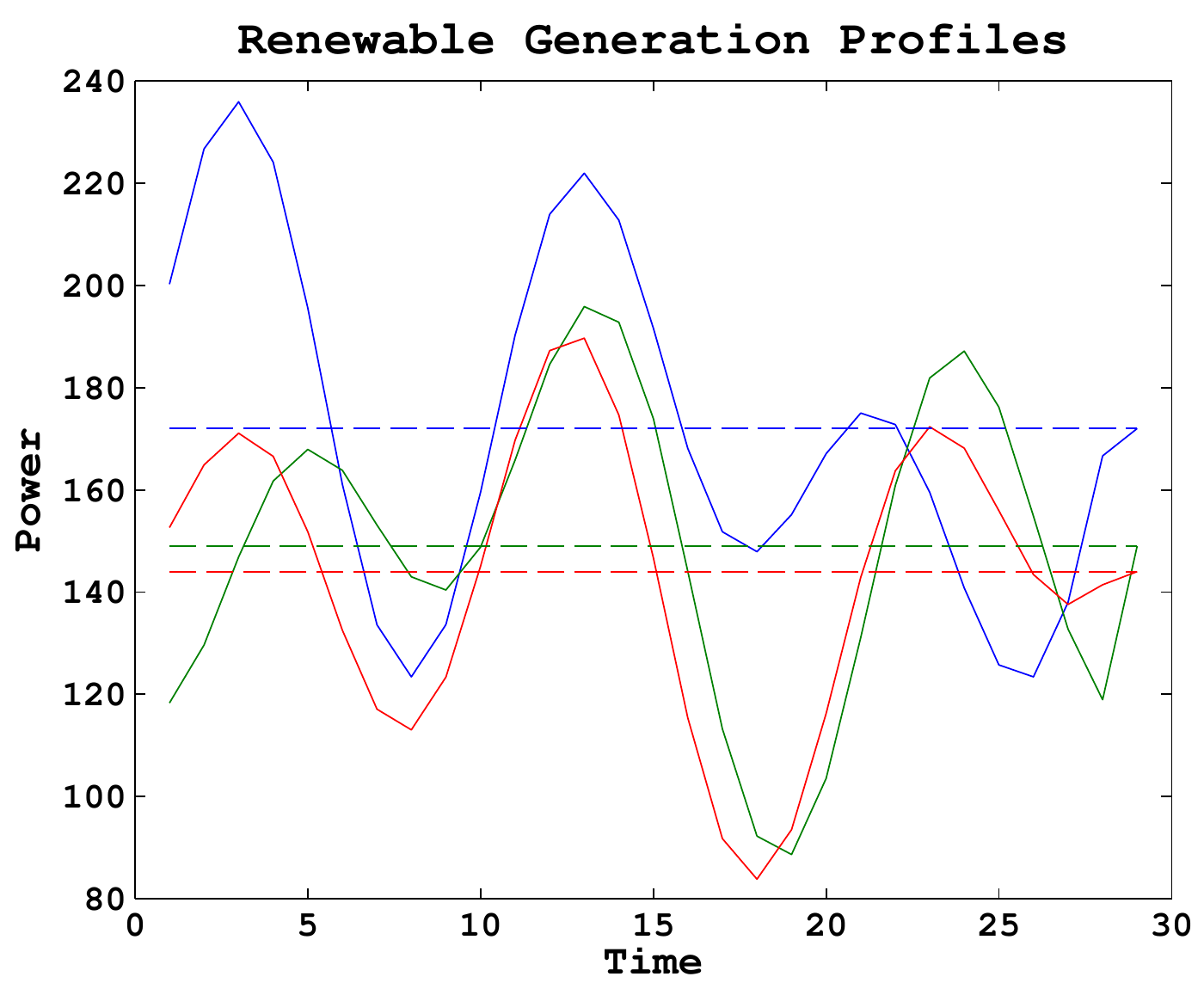}
\caption{Three typical renewable generation profiles used in our study.  The dashed lines are the time averages of the profiles.  } \label{fig:Profiles}
\end{center}
\end{figure}
We seek to develop heuristics to decide how to place storage and size its energy and power capacity. The high-level pseudocode given in algorithm~\ref{alg:OptPlace}, which uses the optimal control algorithm given above, is our initial attempt at this heuristic.  Let $\So$ denote the set of nodes with non-zero storage and $\mathrm{perf}(\So)$ denote the inverse of the total storage power capcity used by the optimal control with access to unlimited storage only at nodes in $\So$.
\begin{algorithm}
\caption{Optimal Sizing and Placement Algorithm}
\label{alg:OptPlace}
\begin{algorithmic}
  \STATE Choose thresholds $\epsilon,\epsilon'$
  \STATE $\So \gets 1:n$
 \REPEAT{}
 \FOR{$k = 1 \to N$}
  \STATE Generate Random Time Series Profiles for the Renewables
  \STATE Solve \eqref{eq:OptProb} with storage only at nodes $\So$ for the given profiles to get optimal dispatch $\popt(t)$
  \STATE $A^k \gets \max_t \abs{\popt(t)}$
 \ENDFOR
  \STATE $A \gets \frac{1}{N} \sum_{k=1}^N A^k$
  \STATE $\gamma \gets \max\{\gamma:\{\mathrm{perf}(\{i \in \So:A_i\geq\gamma\max(A)\})<\mathrm{perf}(\So)+\epsilon\}\} $
  \STATE $\So \gets \{i \in \So:A_i\geq\gamma\max(A)\}$
 \UNTIL{$\gamma\leq\epsilon'$}
\end{algorithmic}
\end{algorithm}
At every iteration of the outer loop, the set of active storage nodes $\So$ is chosen. For each inner loop iteration, we create a trial where a random profile for renewable generation is created at each renewable node (details of how this is done are in Section~\ref{sec:exp} and typical examples of three profiles are shown in Fig.~\ref{fig:Profiles}). For each trial of the renewable profile, the control is allowed access to an unlimited amount of storage power and energy at all the nodes in $\So$ (and $0$ storage at nodes not in $\So$), and we solve the control problem for the optimal storage power dispatch $\popt_{\So}(t)$.  For each trial $k$ and each node in $\So$, we find $A^k$, the maximum of the absolute value of the storage power at that node, which is a measure of the size of the storage power required at that node by the optimal storage dispatch for a given trial $k$.  The values of $A^k$ are averaged over many trials.

Next, the set of storage nodes $\So$ is reduced by picking only those nodes $i \in \So$ such that $A_i\geq\gamma\max(A)$, i.e. only the nodes that require, on average, a storage power capability within a fraction $\gamma$ of the node with the maximum storage power capability.  In essence, we are using our optimal control scheme to identify the nodes that exercise the greatest control over $\costt{\pv,\pvns,\pvs}$. The parameter $\gamma$ determines the severity of the ``cut''.  As $\gamma$ is increased towards one, fewer and fewer nodes remain in $\So$.  As our controllability is reduced, we might expect the level of constraint violations to increase. However, the violations seem to remain almost constant and we instead observe that at a critical value of $\gamma$, the total amount of storage power on the network, as required by the optimal control, suddenly increases.  We set $\gamma$ just below this threshold thereby selecting the reduced set of nodes that performs nearly as well as allowing the optimal control access to all of nodes.

Using the reduced set of nodes from above as our new starting point, the algorithm is repeated (i.e. we regenerate $A$ and find the new maximal reduction in the number of nodes). This outer loop repeats until we are not able to reduce the number of nodes any farther without observing a significant increase in the total amount of storage power.

\section{Simulations}
\label{sec:exp}
\begin{figure}[htb]
\begin{center}
 \includegraphics[scale=0.6]{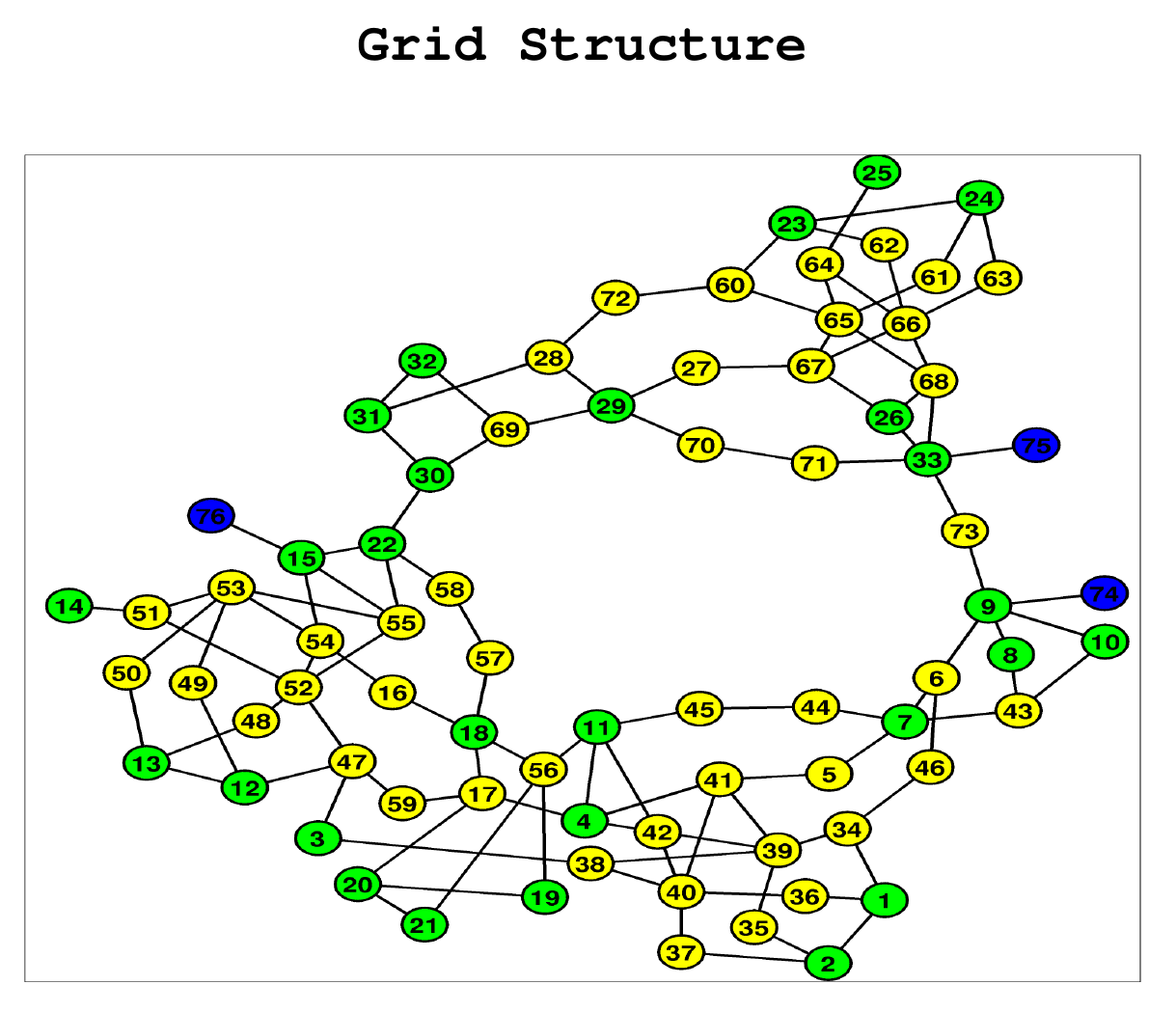}
\caption{Our modified version of RTS-96.  The added renewables are blue, loads are yellow and controllable generators are green.} \label{fig:Grid}
\end{center}
\end{figure}
\begin{figure}[htb]
\begin{center}
 \includegraphics[width=0.45\textwidth,height=.3\textheight]{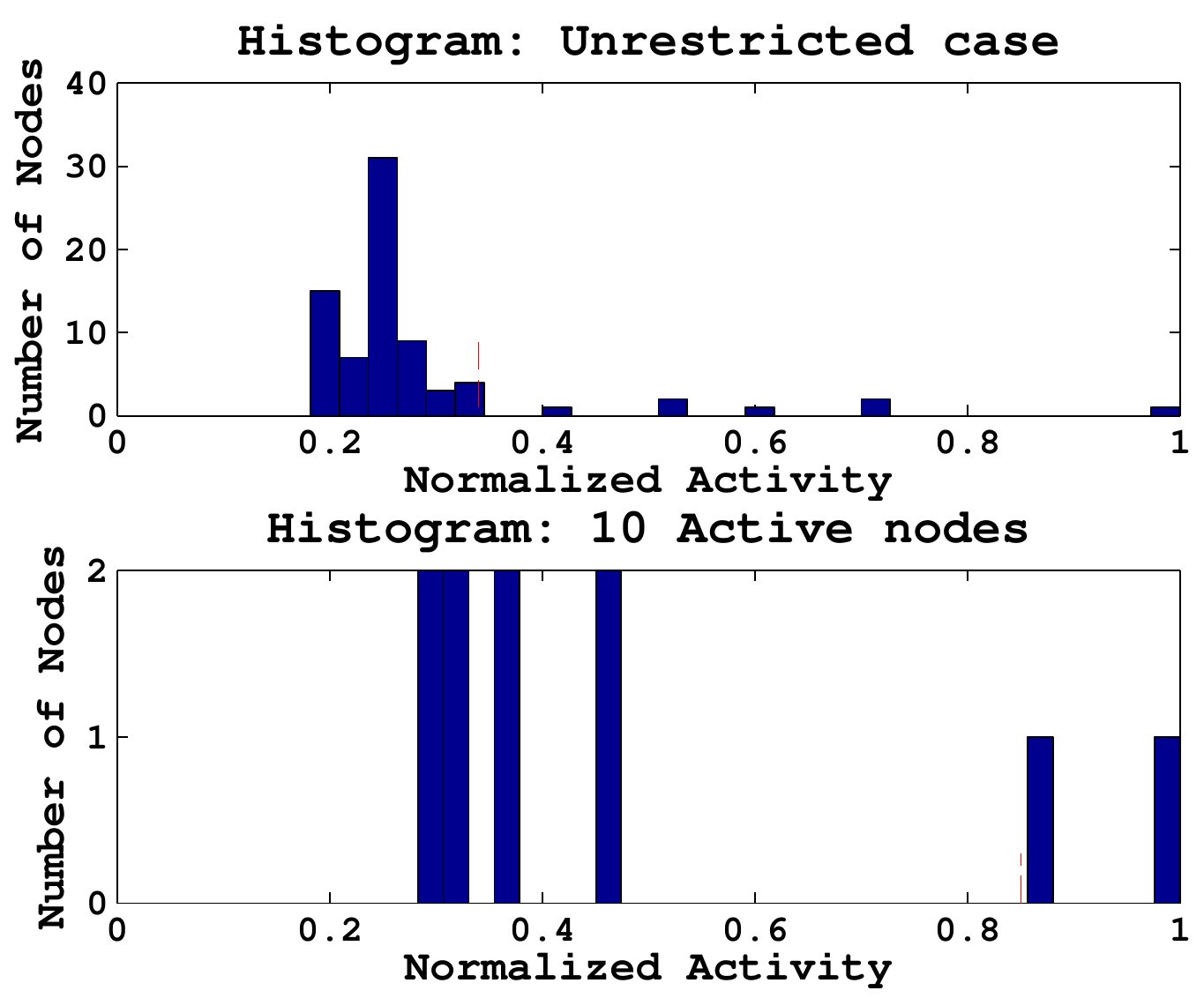}
\caption{Activity Histograms: Red lines mark thresholds used for the reduction in the storage node set.}\label{fig:Hist}
\end{center}
\end{figure}

\begin{figure}[htb]
\begin{center}
 \includegraphics[scale=0.6]{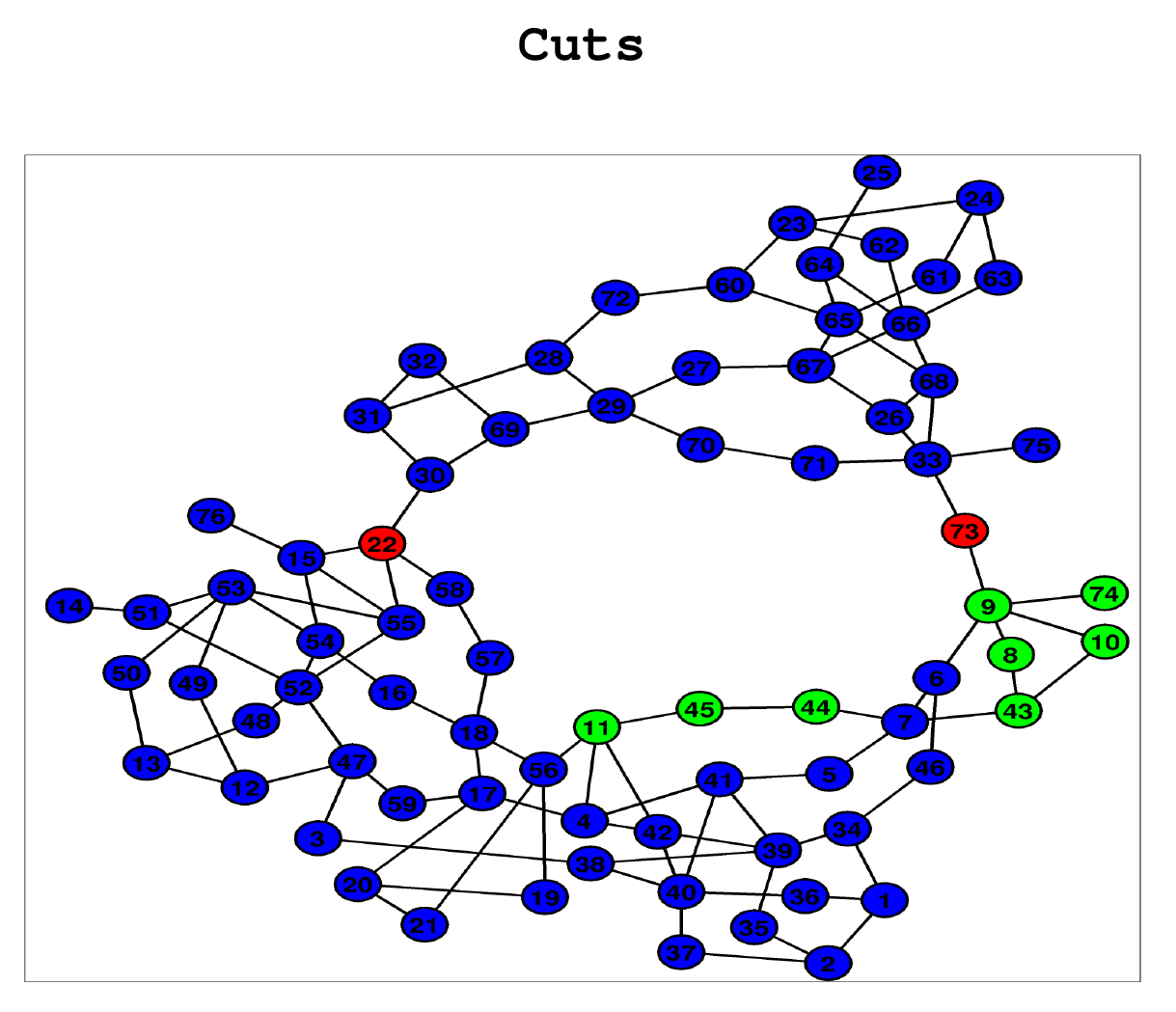}
\caption{Our modified RTS-96 grid showing the three sets of nodes identified by our heuristic.  The minimal set (two nodes) is shown in red.  The ten node set includes the two red nodes and the eight addtional green nodes.  The maximal set includes all of the nodes.} \label{fig:Cuts}
\end{center}
\end{figure}

\begin{figure}[htb]
\begin{center}
 \includegraphics[width=0.4\textwidth]{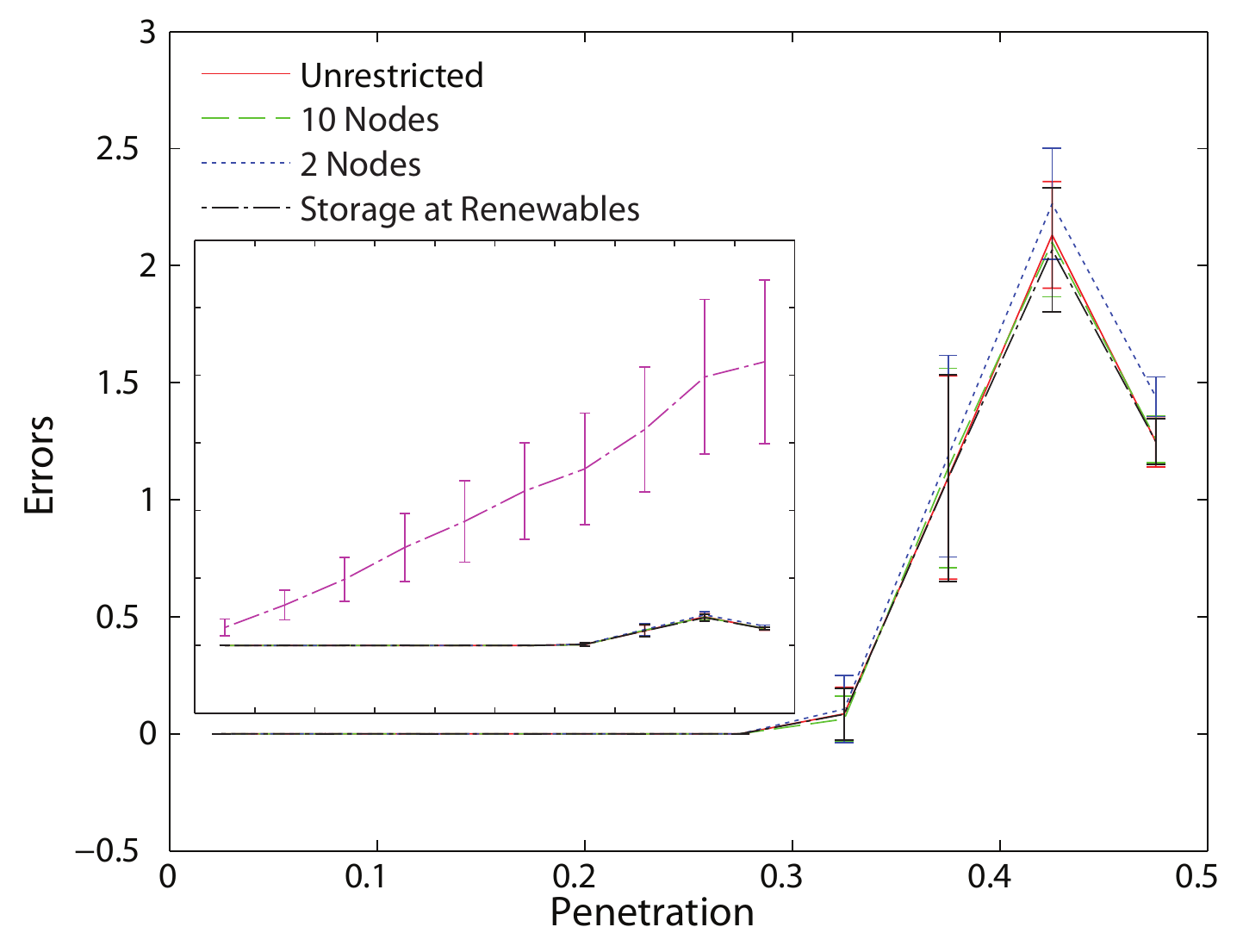}
\caption{The level of constraint violations vs the penetration of renewable generation. } \label{fig:ErrorP}
\end{center}
\end{figure}

\begin{figure}[htb]
\begin{center}
 \includegraphics[width=0.4\textwidth]{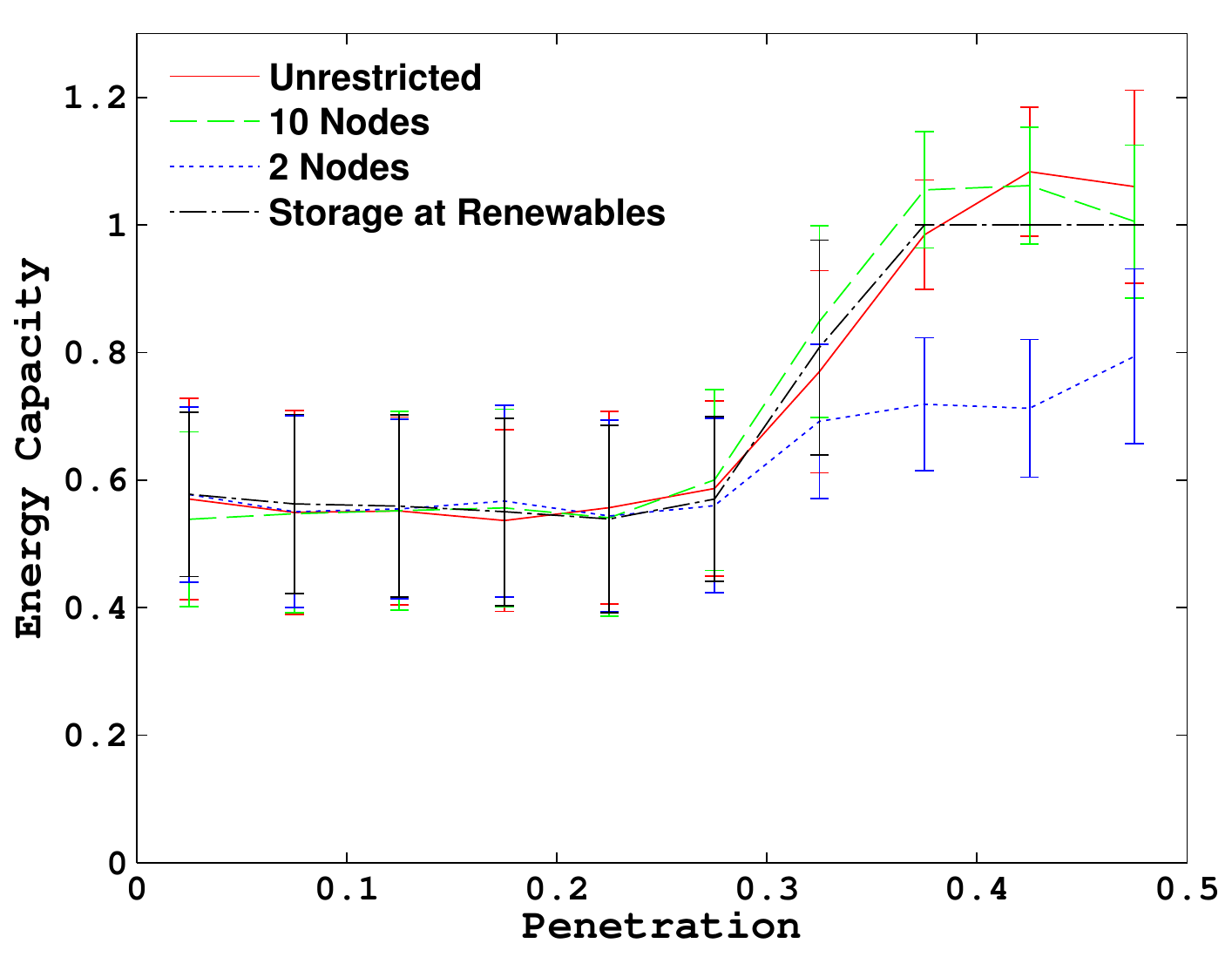}
\caption{The normalized energy capacity of storage in the entire network vs the penetration of renewable generation.} \label{fig:EnergyCapP}
\end{center}
\end{figure}

\begin{figure}[htb]
\begin{center}
 \includegraphics[width=0.4\textwidth]{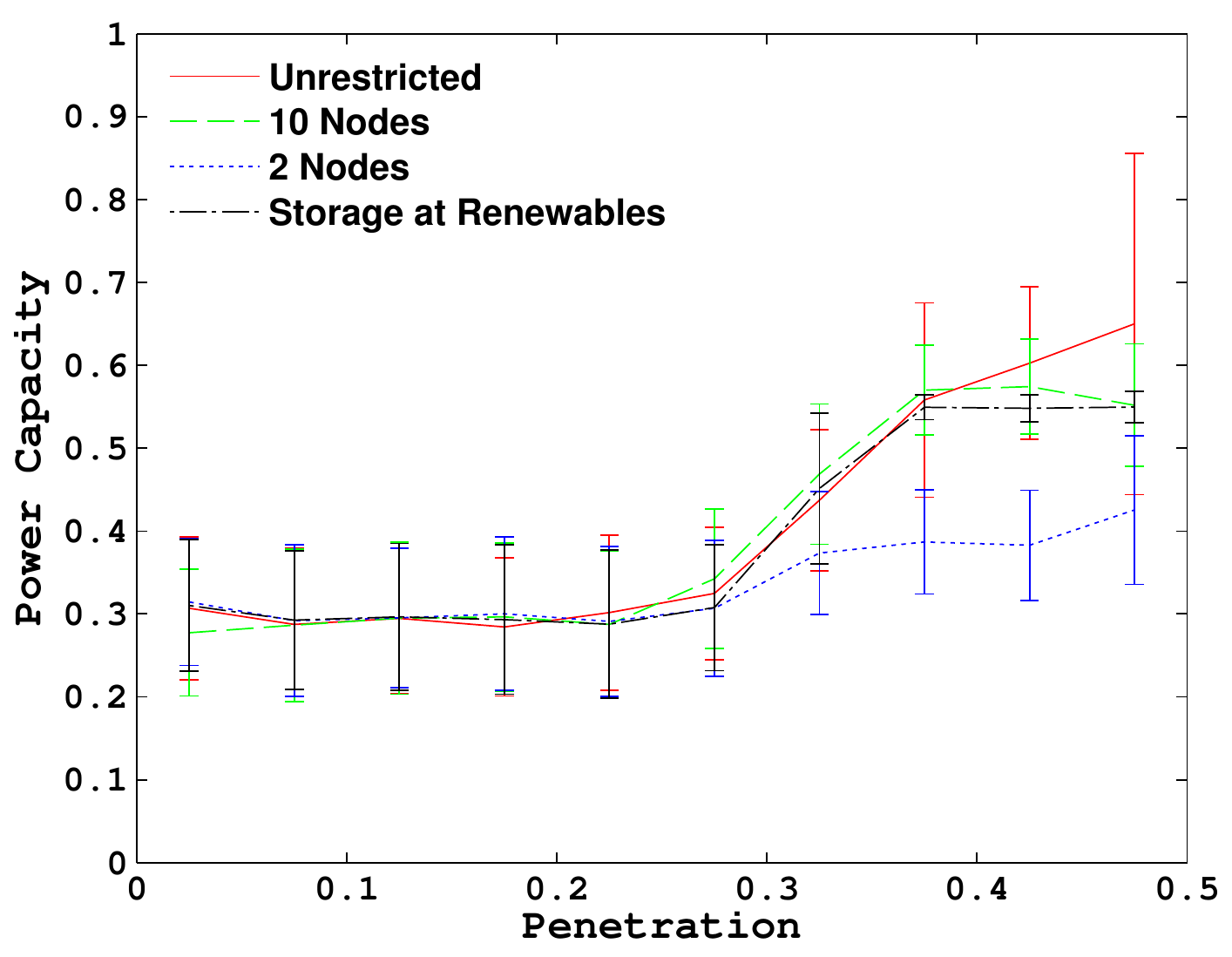}
\caption{The normalized power capacity of storage in the entire network vs the penetration of renewable generation.} \label{fig:PowerCapP}
\end{center}
\end{figure}
% \begin{figure}[htb]
% \begin{center}
%  \includegraphics[scale=0.6]{Figures/Grid.pdf}
% \caption{Grid: Renewables are blue,loads yellow and generators green } \label{fig:Grid}
% \end{center}
% \end{figure}

We tested our optimal control and heuristic for storage placement and sizing on a modified versions of RTS-96\cite{RTS96}.  The grid is shown in Fig.~\ref{fig:Grid}.  Our modification includes the addition of three renewable generation nodes shown in blue.  The capacities of the new lines connecting the renewables to their immediate neighbors are set higher than the capacity of the added renewable generation, otherwise, these lines would be overloaded in nearly every trial.

\subsection{Simulation Setup}
We performed three sets of simulations (corresponding to outer iterations of algorithm \ref{alg:OptPlace}), and for each setting, we perform $N=2000$ trials. We define the renewable penetration value of a trial to be
\[\frac{\sum_{t=0}^{T_f} \sum_{i \in \R} \pvns_i(t)}{\sum_{t=0}^{T_f} \sum_{i \in \Lo} \pvns_i(t)}\].

The penetration is the fraction of the load that is served by the renewable generation. For each trial, a random value uniformly distributed between $0$ and $50\%$ is selected for the penetration level, and the mean values of the renewables are scaled so as to achieve these values. A zero-mean fluctuation is generated around the base value for each renewable:
\[\pvr_i(t)=\sum_k \frac{1}{k} \frac{\int_{0}^\Delta \sin\left(\omega_k(t\Delta+\tau)+\phi_k\right)\,d\tau}{\Delta}\]
where the $\phi_k$ are randomly chosen phases and $\omega_k$ are chosen such that the above sine waves are harmonics of the basic wave with time period $T$.
Since, we discretize the time axis, we use $\frac{\int_{0}^\Delta \sin\left(\omega_k(t\Delta+\tau)+\phi_k\right)\,d\tau}{\Delta}$, the average power in the interval $t\Delta$ to $(t+1)\Delta$.
The final profile is generated by scaling the disturbances so they are comparable to the mean generation (which is typical for turbulent fluctuations of wind) and adding it to the mean: $\pvns_i(t) = \pvnot_i(t)*(1+\pvr_i(t)), i \in \R$. Each instance of the renewable generation profiles we create is statistically independent and uncorrelated with every other instance. We do not intend for a sequence of these instances to approximate an actual time series of wind generation.\\

We start with the unrestricted case (where unlimited storage is available at all nodes) and compute the average activity of each node over all $N=2000$ trials. The histogram of activities is shown in Fig.~\ref{fig:Hist}. The optimal $\gamma$ for the first outer iteration of algorithm \ref{alg:OptPlace} is $0.34$ and the resulting reduced set contains $10$ nodes.  When we tried to remove more nodes from the set, the average level of total network power and energy storage capacity showed a sudden increase indicating that further reductions based on this first stage of planning actually decreased the system performance.

We repeat this process, but this time starting with the 10-node set identified from the first iteration.  Simulating another $N=2000$ trials yields the second histogram in figure \ref{fig:Hist}.  When we restrict our optimal control to just these 10 nodes, their rank order changes relative to the rank order in the unrestricted case.  It is this property that allows our heuristic to make further reductions in the number of nodes without significant impacts to overall performance.

The optimal $\gamma$ for this second iteration is $0.85$ leaving only two nodes for the third interation. The algorithm terminates here as the controllable set cannot be shrunk any further without significance performance degradation. The 3 cuts (all nodes, 10 nodes, 2 nodes) obtained by the algorithm are depicted on the grid in Fig.~\ref{fig:Cuts}.\\

To study how the energy and power capacity requirements of storage depends on the level of renewable penetration, we first must quantify these measures in a consistent fashion.  For a given renewable profile $\pvr(t)$, the optimal storage power dispatch is $\popt(t)$, and $\sopt(t)$ is the resulting energy in storage. We define the energy in the renewable fluctuations to be $\sr(t)=\sum_{\tau=0}^{t-1} \pvr(\tau)\Delta$, i.e. $\sr(t)$ is the energy stored in a (hypothetical) battery that is connected directly to a renewable node and eliminates all fluctuations about the mean renewable generation.   We will only look at changes in energy stored over time, so $\sopt(0)$ and $\sr(0)$ can be arbitrarily chosen.  With these brief prelimiaries, the parameters relevant to our study are:

\begin{itemize}
 \item[]\emph{Normalized Power Capacity}: This quantifies the total power capacity of the storage relative to the sum of maximal power fluctuations over the renewables:
\[\frac{\sum_{j \in \sI} \max_t\abs{\popt_j(t)}}{\sum_{i \in \R} (\max_t \pvr_i(t)-\min_t\pvr_i(t)) }\]
 \item[]\emph{Normalized Energy Capacity}: This quantifies the total energy capacity of the storage relative to the sum of maximal energy fluctuations over the renewables:
\[\frac{\sum_{j \in \sI} (\max_t \sopt_j(t)-\min_t \sopt_j(t))}{\sum_{i \in \R} (\max_t \sr_i(t)-\min_t\sr_i(t)) }\]
 \item[] \emph{Constraint Violations}: This is defined as the sum of all constraint violations (as defined in \eqref{eq:Contraint}) averaged over time.
\end{itemize}

These quantities are all defined for a given renewable configuration. In order to investigate their typical behavior as a function of renewable penetration, we create a 10 equally-spaced intervals between $0$ and $.5$. For each bin $I$, we compute statistics (mean and standard deviation) of all three quantities over trials with a penetration value in the bin $I$. The means an standard deviations for the three quantities are shown as functions of penetration in Figs.~\ref{fig:PowerCapP},\ref{fig:EnergyCapP},\ref{fig:ErrorP}.

\section{Interpretations}
\label{sec:inter}
First, we note that Fig.~\ref{fig:ErrorP} demonstrates the effectiveness of node-specific control for mitigating the impact of renewable generation fluctuations on the network.  The upper curve in the inset in Fig.\ref{fig:ErrorP} shows the level of constraint violations for a system that maintains generation-load balance purely via proportional control, i.e. only the controlled generators responding to a generation-load imbalance in proportion to their capacity.  Without the node-specific control in our model of energy storage, the response of the generators is inflexible, and system begins to experience constraint violations for small values of renewable penetration.  With the node-specific control allowed by our storage model, Figure~\ref{fig:ErrorP} shows that the control can eliminate all constraint violations up to about 30\% renewable penetration for all the node configurations we explored.  We do not attribute this behavior to the storage itself, rather, this appears to be a property of the flexibility of node-specific control.  Although we did not explore this possibility, we believe that if we would have incorporated the controlled generators $\G$ into our control algorithm, we would have seen similar performance up to a renewable penetration of 30\% for a system without any storage.  It is interesting to note that more detailed studies have found similar thresholds for renewable penetration\cite{Meibom2010}.  

Figure~\ref{fig:Cuts} shows the location of the nodes with storage in the reduced-node sets identified by our heuristic in algorithm~\ref{alg:OptPlace}.  The two-node set is perhaps the most illustrative of how optimal placement and sizing of storage in a network does yield itself to simple rules of thumb or intuition.  Intuition may have led to the conclusion that storage should be placed at the site of the fluctuating renewables because it can effectively mitigate the fluctuating power flows without these flows ever being injected into the network.  It seems plausible that this should lead to a minimal amount of storage and quality control, however, the results in Fig.~\ref{fig:Cuts} and \ref{fig:ErrorP} show otherwise.  The nodes in the two-node set do not include either of renewable generation sites.  In fact, even when expanded to the ten-node set, only one of the three renewable sites is included.  These results demonstrate that the nodes that have a high degree of control over congestion caused by renewable fluctuations are not necessarily the renewable nodes themselves.

To further demonstrate this point, we performed an additional simulation where the reduced-note set included just the three renewable sites. Figure~\ref{fig:ErrorP} shows that the quality of control over the constraint violations is the same for this handpicked case as in the two-, ten-, and unrestricted-node sets.  Figures~\ref{fig:EnergyCapP} and \ref{fig:PowerCapP} help to distinguish between these node sets.  For renewable penetrations below about 30\%, the network total power capacity and energy capacity of the storage is nearly the same for each case.  Beyond 30\% where we begin to see a few constraint violations in Fig.~\ref{fig:ErrorP}, the two-node set shows significantly lower total network power and energy capacity compared to the ten or unrestricted-node cases.  The handpicked case is comparable to the ten-node case.  A lower network total capacity implies that the nodes identified in the two-node case have the highest degree of controllability over the network congestion caused by the fluctuating renewables.

A possible reason for the two-node set showing better controllability is evident from Fig.~\ref{fig:Cuts}.  These two nodes sit at the end or the middle of crucial transmission lines that link two major regions in RTS-96.  In reference to Fig.~\ref{fig:Cuts}, we name these regions ``upper'' and ``lower''.  We conjecture that by controlling the power injections at these nodes, our control scheme can effectively control to which region the renewable generation is directed.  Therefore, our control can direct the renewable fluctuations to the region that, at that point in time, has the greatest ability to assist in mitigation.  In addition to simply absorbing the fluctuations locally, the storage may be used to redirect the renewable fluctuations throughout the network to leverage other resources to assist in the control.
\section{Conclusions and Future Work}
\label{sec:con}
In this initial study, we have developed and demonstrated a method for coupling operations with system expansion planning for the optimal placement and sizing of storage in a grid with a significant penetration of time-intermittent renewable generation.  Operations are incorporated into the planning via simulations of an optimal control scheme that uses perfect renewable generation forecasts to dispatch energy storage to eliminate violation of network and generation constraints.  By simulating the system on a relatively short time scale, we build up significant statistics on the storage activity at all network nodes.  We have developed a heuristic that uses these statistics to reduce the number of storage-enabled network nodes while maintaining or improving system performance.

Somewhat unexpectedly, our method does not preferentially place energy storage at the nodes with renewable generation.  Instead, our method apparently favors nodes at critical junctions between major subcomponents of the network.  We have proposed that these nodes provide for enhanced controllability because, in addition to simply buffering the fluctuations of the renewables, controlled power injections at these nodes can modify overall network flows and direct fluctuating power flows to regions that are better positioned to mitigate them.  

There is much follow on work needed to expand the concept presented in this manuscript and to verify some of its conjectures.  A few examples include:
\begin{itemize}
\item Additional networks should be considered, including different configurations of renewable generation on the present network, to investigate whether our node-reduction heuristic is robust.
\item The explanation of the placement of storage proposed in this work should be verified by correlating the response of the storage with nearby power flows.
\item  The operational simulation should be made more realistic by incorporating ramping constraints on the controllable generation.  
\item  The method should be extended to consider alternatives to proportional control for changing the outputs of the controllable generators.
\end{itemize} 

% %{\color{red} !!!!! Acknowledgments are needed only for the final submission - hide before submitting}
\section{Acknowledgments}

We are thankful to the participants of the "Optimization and
Control for Smart Grids" LARD DR project at Los Alamos
and Smart Grid Seminar Series at CNLS/LANL for multiple
fruitful discussions. The work at LANL was carried out under the auspices of the National
Nuclear Security Administration of the U.S. Department of Energy at Los
Alamos National Laboratory under Contract No. DE-AC52-06NA25396.

{\small
\bibliographystyle{unsrt}
\bibliography{SmartGrid}
}

\end{document}